\documentclass[conference]{IEEEtran}
\IEEEoverridecommandlockouts
\usepackage{cite}
\usepackage{amsmath,amssymb,amsfonts}
\usepackage{algorithmic}
\usepackage{graphicx}
\usepackage{textcomp}
\usepackage{colortbl}
\usepackage{fancyhdr}
\usepackage{booktabs}
\usepackage[keeplastbox]{flushend}
\usepackage{marvosym}
\ifCLASSOPTIONcomsoc
\usepackage[caption=false, font=normalsize,labelfont=sf,textfont=sf]{subfig}
\ELSE
\usepackage[caption=false,font=footnotesize]{subfig}
\fi
\usepackage[usenames,dvipsnames,svgnames]{xcolor}

\def\BibTeX{{\rm B\kern-.05em{\sc i\kern-.025em b}\kern-.08em
    T\kern-.1667em\lower.7ex\hbox{E}\kern-.125emX}}
 
\DeclareMathOperator{\diff}{d}    
\newcommand{\ddt}{\frac{\diff}{\diff t}}

\usepackage{tikz}
\usepackage{circuitikz}
\usepackage{pgfplots}
\usetikzlibrary{arrows}
    
\begin{document}

\title{Dispatchable Virtual Oscillator Control for Decentralized Inverter-dominated Power Systems: Analysis and Experiments\\
%{\footnotesize \textsuperscript{*}Note: Sub-titles are not captured in Xplore and should not be used}
\thanks{This work was authored in part by the National Renewable Energy Laboratory, operated by Alliance for Sustainable Energy, LLC, for the U.S. Department of Energy (DOE) under Contract No. DE-AC36-08GO28308. Funding was provided in part by the DOE Office of Energy Efficiency and Renewable Energy Solar Energy Technologies Office grant number DE-EE0000-1583. This project has received funding in part from the European Union's Horizon 2020 research and innovation program under grant agreement No$^\circ$ 691800. {This paper reflects only the authors' views.} The views expressed in the article do not necessarily represent the views of the DOE{, the U.S. Government, or the European Commission.  The European Commission is not responsible for any use that may be made of the information it contains.}
}
}

\author{\IEEEauthorblockN{Gab-Su Seo\IEEEauthorrefmark{1}, Marcello Colombino\IEEEauthorrefmark{1}, Irina Subotic\IEEEauthorrefmark{2}, Brian Johnson\IEEEauthorrefmark{3}, Dominic Gro\ss\IEEEauthorrefmark{2}, and Florian D{\"o}rfler\IEEEauthorrefmark{2}}
\IEEEauthorblockA{\IEEEauthorrefmark{1}Power Systems Engineering Center, National Renewable Energy Laboratory, Golden, CO 80401, USA\\
Email: Gabsu.Seo@nrel.gov, Marcello.Colombino@nrel.gov}
\IEEEauthorblockA{\IEEEauthorrefmark{2}Automatic Control Laboratory, ETH Zurich, Zurich 8092, Switzerland\\
Email: subotici@ethz.ch, grodo@ethz.ch, dorfler@ethz.ch}
\IEEEauthorblockA{\IEEEauthorrefmark{3}Department of Electrical and Computer Engineering, University of Washington, Seattle, WA 98195, USA\\
Email: brianbj@uw.edu}
}

\maketitle

\begin{abstract}
This paper presents an analysis and experimental validation of dispatchable virtual oscillator control (dVOC) for inverter-dominated power systems. dVOC is a promising decentralized control strategy that requires only local measurements to induce grid-forming behavior with programmable droop characteristics. It is dispatchable--i.e., the inverters can vary their power generation via user-defined power set-points and guarantees strong stability. To verify its feasibility, a testbed comprising multiple dVOC-programmed inverters with transmission line impedances is designed. With an embedded synchronization strategy, the dVOC inverters are capable of dynamic synchronization, black start operation, and transient grid voltage regulation with dynamic load sharing, and real-time-programmable droop characteristics for backward compatibility. All these features are experimentally verified.
\end{abstract}


\begin{figure*}[t!]
\centerline{\includegraphics[width = 2.0\columnwidth]{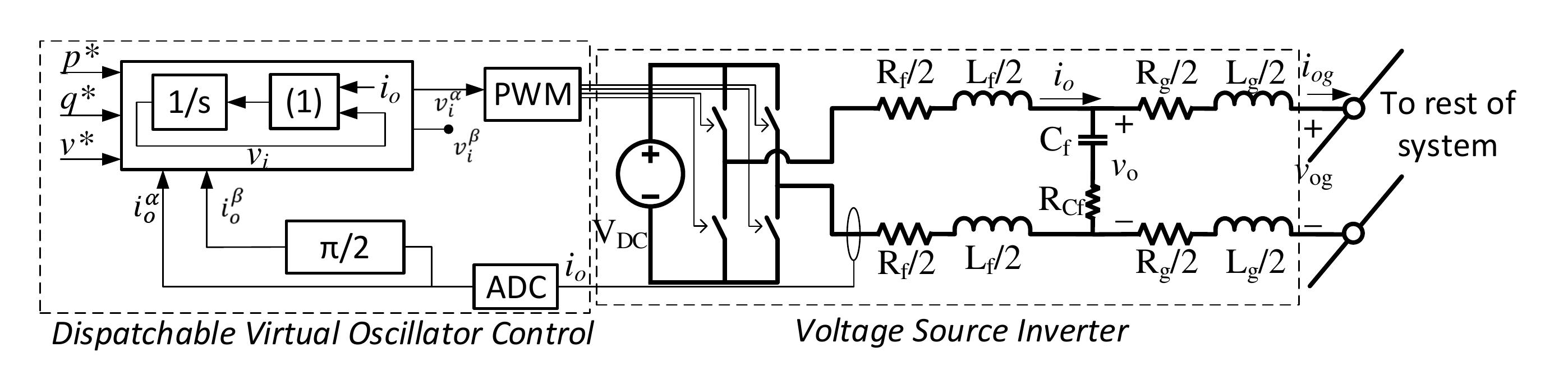}}
\caption{Schematic of dVOC inverter for decentralized inverter-dominant grid}
\label{fig:dvoc:schematics}
\end{figure*}

\section{Introduction}
The electric power grid is undergoing exceptional changes with increasing penetration of inverter-based renewable generation~\cite{iea2016, nrel_renewable_book_2016}. As of now, stability and system-wide synchronization of the grid is achieved with traditional synchronous generators and their controls. Conventionally, power inverters are controlled in a ``grid following" fashion. This means that they are programmed to estimate the (already stable) grid frequency and regulate their injected current to track pre-determined power set-points. In a purely inverter-based grid, without synchronous generators, these control strategies are no longer suitable as they cannot maintain stability and synchronization. To obviate this problem and envision a stable inverter-based grid, “grid-forming” control methods have been proposed~\cite{torres2015synchronization, arghir2018grid, liu2016comparison}. A grid-forming inverter is not limited to power tracking but acts as a controlled voltage source that can change its power output (thanks to storage or curtailment), and is controlled to contribute to the stability of the grid. It is highly desirable that inverter control strategies are decentralized, i.e., they rely on local measurements only~\cite{kroposki2017achieving} as this allows for plug-and play capabilities. We envision that, with the advent of renewable generation and grid-forming control, portions of the grid can be islanded if needed without loosing synchronization and load sharing. This improves modularity and thus resiliance of the grid to natural disasters or cyber-physical attacks.

Most of the common approaches of grid-forming control focus on droop control~\cite{MCC-DMD-RA:93,Guerrero2015,simpson2017voltage}. Its simple implementation and backward compatibility~\cite{chandorkar1993control}, make droop control a desirable solution accepted by utilities and practitioners. {However, the associated phasor models are well-defined only near the synchronous steady-state \cite{johnson2017comparison}.} Other popular approaches are based on mimicking the physical characteristics and controls of synchronous machines~\cite{zhong2011synchronverters,jouini2016grid,SDA-SJA:13}. While strategies based on machine-emulation and droop are compatible with legacy power systems, they use an energy conversion interface (an inverter) with fast actuation but almost no inherent energy storage to mimic another interface (a generator) with slow actuation but significant energy storage (in the form of rotational inertia). It is unclear if this is a viable option, especially taking into account the {limited energy storage} and tight over-current constraints of power inverters~\cite{ENTSOE16}. 

Recently proposed solutions based on virtual oscillator control (VOC) feature enhanced dynamic performance and maintain an embedded droop control law {\cite{MS-FD-BJ-SD:14b}} close to steady-state. This implies a superior voltage regulation performance with respect to standard droop control, while keeping load sharing capabilities~\cite{johnson2014synchronization, johnson2017comparison}. In addition to these benefits, ease of synchronization allows VOC to be a promising candidate for microgrids; however, it is unclear how to “dispatch” VOC inverters--i.e., reconfigure the inverters' power injections as they have no programmable power set-points. 
This paper investigates a recently proposed dispatchable virtual oscillator control (dVOC) with the following desirable features: 
\begin{itemize}
    \item [i)]  dVOC is dispatchable, i.e. it allows for the user to specify power set-points for each inverter. 
    \item [ii)] Given no set-points, dVOC subsumes VOC control and therefore it inherits all its favorable dynamical properties.
    \item [iii)] Under the assumption that the set-points are consistent with the AC power flow equations and other technical assumptions, dVOC renders the grid globally asymptotically stable with respect to the desired solution of the AC power-flow~\cite{colombino2017global}.
\end{itemize}
The contribution of this paper is twofold: first, we analyze the inherent droop characteristics that dVOC shows around a synchronous trajectory. This allows practitioners who are accustomed to working with droop controllers to understand the effect of the dVOC parameters on its steady-state behavior. Finally, we present an experimental validation of dVOC's feasibility for an inverter-based microgrid, where we validate its stability, load sharing, and droop properties in a two-inverter system. 

The paper is structured as follows. In Section~\ref{Sec:ii} we introduce the dVOC control; in Section~\ref{Sec:iii}, we derive the dVOC droop characteristics; in Section~\ref{Sec:iv}, we present the experimental setup and preliminary validation results for dVOC, and in Section~\ref{Sec:v} we conclude by summarizing the results and outlining the planned future work.

% %%%%%%%%%%%%%%%%%%%%%%% FIGURE %%%%%%%%%%%%%%%%%%%%%%%%%%%%%%%%
% \definecolor{mycolor1}{rgb}{0.00000,0.44700,0.74100}%
% \definecolor{mycolor2}{rgb}{0.85000,0.32500,0.09800}%
% \definecolor{mycolor3}{rgb}{0.92900,0.69400,0.22500}%

% \begin{figure}[htbp]
% \begin{center}
% \begin{circuitikz}[american voltages]

% \ctikzset{bipoles/resistor/height=0.15}
% \ctikzset{bipoles/resistor/width=0.4}

% \ctikzset{bipoles/generic/height=0.15}
% \ctikzset{bipoles/generic/width=0.4}

% \ctikzset{bipoles/length=.6cm}

% \coordinate (I1) at (-2,0);
% \coordinate (I2) at (2,0);
% \coordinate (I3) at (0,0);

% \draw ($ (I1) + (0,-1) $) node [ground] (g1) {};
% \draw (g1) to[sV, l=$v_1$, color=mycolor1] (I1);

% \draw ($ (I2) + (0,-1) $) node [ground] (g2) {};
% \draw (g2) to[sV, l=$v_2$, color=mycolor2] (I2);

% \draw ($ (I3) + (0,-1) $) node [ground] (g3) {};
% \draw (g3) to[R, l=$R_{l}$] (I3);

% \node (I13) at ($(I1)!0.5!(I3)$)[label={[xshift=0cm, yshift=-0.2cm]$Z_{1l}$}]{};
% \draw 		(I1) 
% to[R,*-] (I13)
% to [L,-*] (I3);

% \node (I32) at ($(I3)!0.5!(I2)  $) [label={[xshift=0cm, yshift=-0.2cm]$Z_{2l}$}]{};
% \draw 		(I3) 
% to[R,*-] (I32)
% to [L,-*] (I2);

% \end{circuitikz}
% \caption{Inverter based grid}
% \label{fig.inverter.example}
% \end{center}
% \end{figure}
% %%%%%%%%%%%%%%%%%%%%%%%%%%%%%%%%%%%%%%%%%%%%%%%%%%%%%%%%%%%%%%%%%%%%%%%%%%%%%%%%%%%%%

\section{{dVOC} for Inverters}\label{Sec:ii}
This section reviews the basics of dVOC. A detailed description of the control strategy can be found in~\cite{colombino2017global, gross2018effect}.

\subsection{Dispatchable Virtual Oscillator Control}

% Fig. 1 shows an inverter dominated grid, where node i consists of a voltage source $v_{o,i}$ essentially made of an inverter (multiple inverters can be aggregated as a single representative~\cite{purba2017reduced}) and/or a local load, $R_{l,i}$. The nodes i and j are connected with the transmission lines, $Z_{ij}$. 
% In this work, we consider an inverter-dominant grid with resistive loads and inductive-resitive transmission lines. In an inverter-dominant grid, power inverters are not only responsible for supplying power to the loads, but also for stability and synchronization.

Dispatchable VOC is a dencentralized ``grid forming" control strategy designed to achieve synchronization of an inverter-dominant grid, while maintaining a level of control on the power injections and voltage level of each inverter. When applying dVOC, each inverter monitors its output current $i_o$ and, using a PWM strategy, regulates its terminal voltage vector $v_i=[v_i^\alpha,\, v_i^\beta ]^\top$, (in the $\alpha-\beta$ coordinate frame) to follow the dVOC control law:
\begin{equation}\label{eq.dvoc}
\ddt v_i = \omega_0 J v_i+\eta \left(K_i v_i-R(\kappa) i_{o,i}+\alpha\phi_i (v_i ) v_i \right),
\end{equation}
where $i_{o,i}=[i_{o,i}^\alpha,\,i_{o,i}^\beta ]^\top$ is the measurement of the inverter current (for single phase signals the $\beta$ component is reconstructed using a Hilbert transform), $\omega_0$ is the nominal grid frequency, the matrix 
\[
R(\kappa):=\begin{bmatrix}\cos(\kappa)& -\sin (\kappa)\\ \sin(\kappa)& \cos(\kappa)\end{bmatrix}
\]
 is a 2D rotation matrix, ${J:=R(\pi/2)}$,  
 \[
 K_i := \frac{1}{v_i^{\star 2} } R(\kappa) \begin{bmatrix} p_i^\star& q_i^\star\\ -q_i^\star& p_i^\star \end{bmatrix}, \quad \phi_i (v_i ):=\frac{v_i^{\star2}-\|v_i \|^2}{v_i^{\star2}},
 \]
 the operator $\|\cdot\|$ is the Euclidean norm, the quantities $\eta > 0$, $\alpha > 0$ and $0\le\kappa\le\pi$ are the design parameters, and $p_k^\star,q_k^\star$, and $v_k^\star$ are the active power, reactive power, and voltage magnitude {set-points}, respectively. {The parameter $\kappa$ can be used to adjust the controller to adapt to the line parameters, $\kappa=0$ corresponds to resistive lines and $\kappa=\pi/2$ corresponds to inductive lines.} Fig.~\ref{fig:dvoc:schematics} shows a schematic of an inverter {implementing} dVOC.
 
 \subsection{Interpretation of the dVOC Controller}
 Equation~\eqref{eq.dvoc} represents the dVOC controller introduced in~\cite{colombino2017global}. In order to give an intuitive interpretation to dVOC, we can write~\eqref{eq.dvoc} as 
\begin{align}\label{eq.dvoc2}
\ddt v_i = \omega_0 J v_i {+ \eta e_{\theta,i}(v_i,i_{o,i}) + \eta \alpha e_{v,i}(v_i)},
\end{align}
where $\omega_0 J v_i$ is the standard equation of a harmonic oscillator in rectangular coordinates, {$e_{\theta,i}(v_i,i_{o,i}) = K_i v_i - R(\kappa)i_{o,i}$ represents a ``phase error" term (see \cite[Sec. II. D.]{gross2018effect}), and $e_{v,i}(v_i) = \phi_i(v_i)v_i$ represents a magnitude error term. Note, that the magnitude error vanishes when $\|v_i\|= v_i^\star$. Moreover, the ``phase error" $e_{\theta,i}(v_i,i_{o,i})$  vanishes when the voltage levels and power injections of the inverters corresponds to the set-points. To see this note that}
\begin{align*}
\frac{1}{v_i^{\star 2} }  \begin{bmatrix} p_i^\star& q_i^\star\\ -q_i^\star& p_i^\star \end{bmatrix} v_i - i_{0,i} = 0   
\end{align*}
whenever $v_i^\top i_{o,i} = p^\star_i$ and $v_i^\top J i_{o,i} = q^\star_i$. In~\cite{colombino2018global,gross2018effect}, it is shown that, if all inverters in the grid implement (1), and all set-points are consistent with the AC power-flow equations, and further technical assumptions are satisfied, the inverter-based grid is (almost) globally asymptotically stable with respect to the desired power-flows. This means that regardless of the initial conditions, the inverters synchronize and reach the desired set-points. Furthermore, when the set-points are inconsistent with the power-flow equations, dVOC presents droop-like characteristics that enable {grid} synchronization while also giving trade-offs between power imbalance and reactive power versus frequency and voltage. These droop characteristics will be discussed in depth in the next section.

\begin{figure}[t!]
\centerline{\includegraphics[width = 1.0\columnwidth]{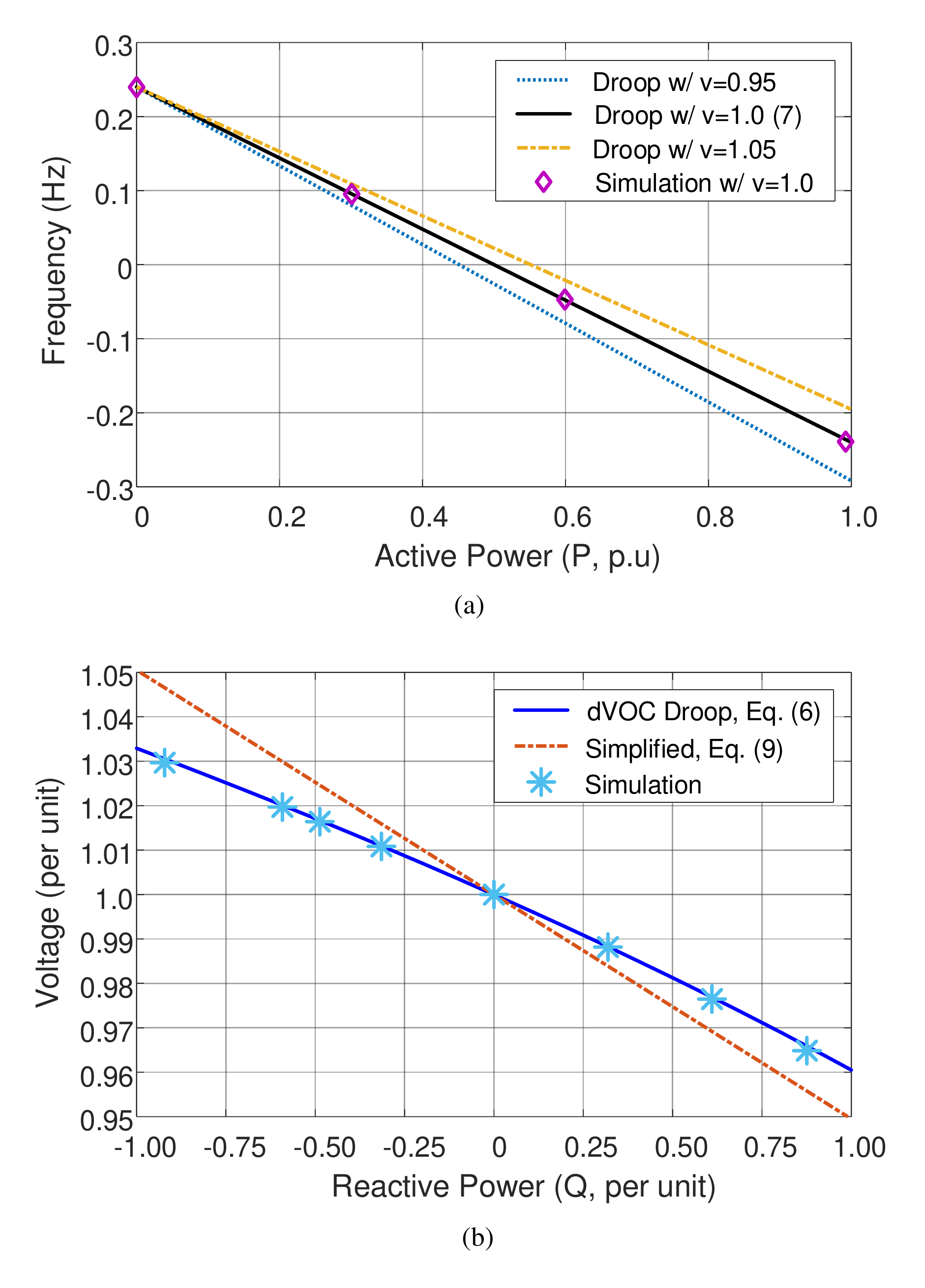}}
\caption{Droop characteristics of dVOC: (a) frequency to active power droop, (b) voltage to reactive power droop. We observe that the {theoretical} droop curves fit very well with high fidelity simulations. Experimental validation of the droop curves is part of our current and future efforts.}
\label{fig:droopcurves}
\end{figure}

\section{Droop Characteristics of {dVOC}} \label{Sec:iii}
Droop control is the most popular decentralized grid-forming control strategy as it is backward compatible with generators and enables power sharing~\cite{johnson2016synthesizing}. 
%With conventional droop law, an inverter’s frequency and voltage magnitude are determined by active and reactive power mismatches as follows
{With conventional droop control, the frequency and magnitude of the voltage of an inverter are determined by active and reactive power mismatches as follows}
\begin{align}
\ddt\theta &=\omega_i=\omega_0+k_p (p_i^\star-p_i),\\ 					
\ddt \|v_i\| &=-\|v_i\|+v_i^\star+k_q (q_i^\star-q_i), 					
\end{align}
where $k_p$ and $k_q$ are droop coefficients. If we consider~\eqref{eq.dvoc} in polar coordinates~\cite{gross2018effect}, we observe that dVOC presents the following nonlinear droop behavior
\begin{align}\label{eq.nonlinear.droop}
\begin{split}
\ddt
\begin{bmatrix}
\|v_i\|\\
\theta_i
\end{bmatrix} &= \eta 
\begin{bmatrix}
\|v_i\| & 0 \\
0 & 1 
\end{bmatrix}
R(\kappa)
\begin{bmatrix}
\frac{p_i^\star}{v_i^{\star2}} - \frac{p_i}{\|v_i\|^2}\\
-\left(\frac{q_i^\star}{v_i^{\star2}} - \frac{q_i}{\|v_i\|^2}\right)
\end{bmatrix}\\
&+
\begin{bmatrix}
\frac{\eta \alpha}{v_i^{\star2}} \left( v_i^{\star2} -\|v_i\|^2  \right)\|v_i\| \\
\omega_0
\end{bmatrix}.
\end{split}
\end{align}
Choosing $\kappa = \pi/2$ in~\eqref{eq.nonlinear.droop} yields the final dVOC droop $(V-Q,\, \omega-P)$
\begin{align}\label{eq.final.droop}
\begin{split}
\ddt\theta_i  &= \omega_0+\eta\left(\frac{p_i^\star}{v_i^{\star2}} -\frac{p_i}{\|v_i\|^2} \right)\\
\ddt \|v_i\| & =  \eta \!\left(\frac{q_i^\star}{v_i^{\star2}} -\frac{q_i}{\|v_i\|^2} \right) \|v_i\| \\
& + \frac{\eta \alpha}{v_i^{\star2}} \left( v_i^{\star2} -\|v_i\|^2  \right)\!\|v_i\|.
\end{split}
\end{align}
On the contrary, $\kappa = 0$ yields the {resistive} droop characteristics $(V-P, \omega-Q)$ which is similar to that of VOC~\cite{johnson2017comparison}. It is noteworthy that the system is provably stable when $\kappa = \tan^{-1}\left(\frac{\omega_0L}{R}\right)$, where the inductance/resistance ratio is assumed constant across all transmission lines~\cite{gross2018effect}. We notice both in simulation and in our experimental result that the system remains stable even when this assumption does not hold.
By approximating  $\|v_i\| \approx v_i^\star$, i.e., %small voltage deviation,
{assuming a small voltage magnitude deviation}, a more intuitive droop characteristic can be derived from~\eqref{eq.final.droop} as 
 \begin{align}
\ddt\theta_i =\omega_i & \approx \omega_0+\frac{\eta}{v_i^{\star 2}} (p_i^\star-p_i),\label{eq.droop.p.approx} \\ 			
\ddt \|  v_i\|  &\approx\frac{\eta}{v^\star}(q_i^\star-q_i)+\eta \alpha (v_i^\star-\|  v_i\| ). \label{eq.droop.q.approx}			
\end{align}
Note that at steady-state,~\eqref{eq.droop.q.approx} becomes 
\begin{align}
  \|  v_i\| \approx v_i^\star + \frac{1}{\alpha v_i^\star}(q^\star-q). \label{eq.droop.v.approx}
\end{align}
To illustrates the droop relationships of dVOC, Fig. {\ref{fig:droopcurves}} displays an example design with $p^\star=0.5$ p.u., $q^\star=0$ p.u., $\eta=43.43$  $\Omega$rad/sec and $\alpha=0.9722$ $\mho$  and $\kappa=\pi/2$. As expected from~\eqref{eq.final.droop}, dVOC droop features voltage dependent $\omega-P$ droop, which is different from conventional droop, and non-linear $V-Q$ droop with reduced voltage droop by $q$ change. {We emphasize that dVOC has a local droop characteristic while guaranteeing that network of inverters controlled by dVOC  (almost) globally converges to a pre-defined solution of the AC power flow equations.}

\section{Implementation of dVOC and Experimental Validation}\label{Sec:iv}
\begin{figure}[t!]
\centerline{\includegraphics[width = 1.0\columnwidth]{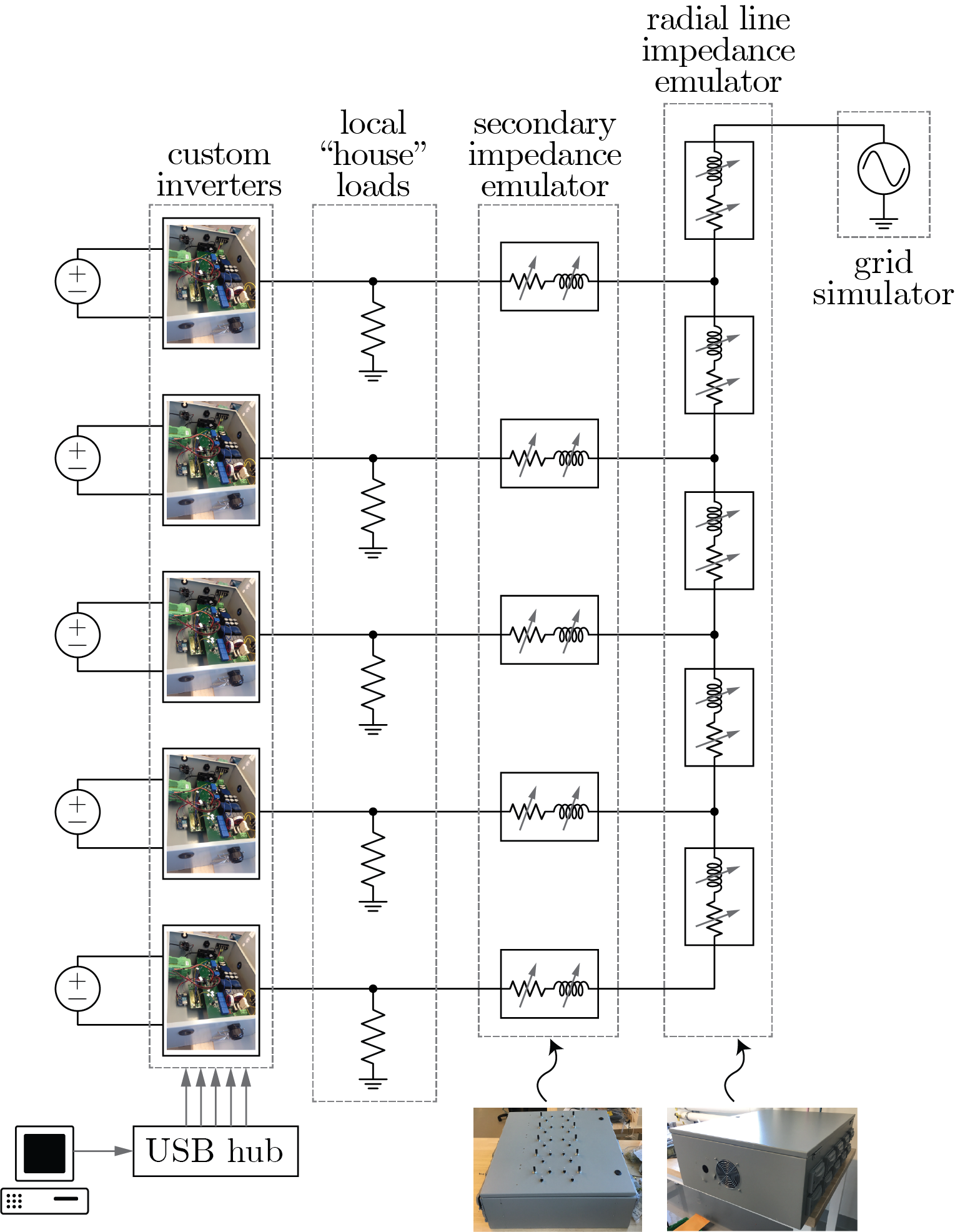}}
\caption{Inverter-dominant grid testbed developed at NREL: Up to five inverters can be connected through a fully  configurable impedance emulator to either a load or a grid simulator. For all experiments presented in this paper we consider two inverters connected in parallel to supply a resistive load. The full potential of the experimental setup will be exploited as part of current and future work.}
\label{fig.testbed}
\end{figure}

\begin{table}
\caption{Design parameters for the experimental validation.}
\label{tab.parameters}
\begin{center}
\begin{tabular}{ @{} cc @{} } 
\bottomrule
\multicolumn{2}{>{\columncolor[gray]{.95}}c}{ dVOC Inverter parameters} \\
\toprule
% & & \multicolumn{3}{c}{Feedback iterations} \\
Item & Design Selection   \\
\toprule 
oscillator param.  &  $\eta=21.71$ $\Omega$ rad/sec, $\alpha=0.9722\, \mho$, $\kappa=\pi/2$\\
set-points &  $p^\star=500$ W, $q^\star= -125$ var, $v^\star=120\, V_{\text{rms}}$\\
controller & 320F28379D, Texas Instruments \\
switching freq. & 32 KHz\\
rated power & 1 VA\\
filter params. & $L_f = 1$ mH, $C_f = 24$ $\mu$F, $L_g = 0.2$ mH\\
\toprule
\bottomrule
\end{tabular}
\end{center}
\end{table}

To verify the concept of dVOC, a hardware testbed shown in Fig.~\ref{fig.testbed} was built. The specifications of the inverter hardware and dVOC parameters are tabulated in Table~\ref{tab.parameters}. The reactive power set-point, $q^\star=-125$ var was set to incorporate the internal reactive power consumption by the filter capacitor $C_f$. As illustrated in Fig~\ref{fig.testbed}, the testbed was constructed to emulate different system condition such as grid connected and islanded conditions with different distribution/transmission line impedances with up to five $1$-kVA inverters. A test scenario was designed to verify collective grid regulation with multiple dVOC inverters. This section provides experimental results for 
\begin{enumerate}
\item [i)] black start operation using a dVOC inverter from dead grid under loaded condition, 
\item [ii)] dynamic synchronization and load sharing of the inverters under loaded condition, 
\item [iii)] load transient performance with two inverters active,
\item [iv)] real-time set-point update (dispatch) operation.
\end{enumerate}

\subsection{Black Start}

\begin{figure}[t!]
\centerline{\includegraphics[width = 1.0\columnwidth]{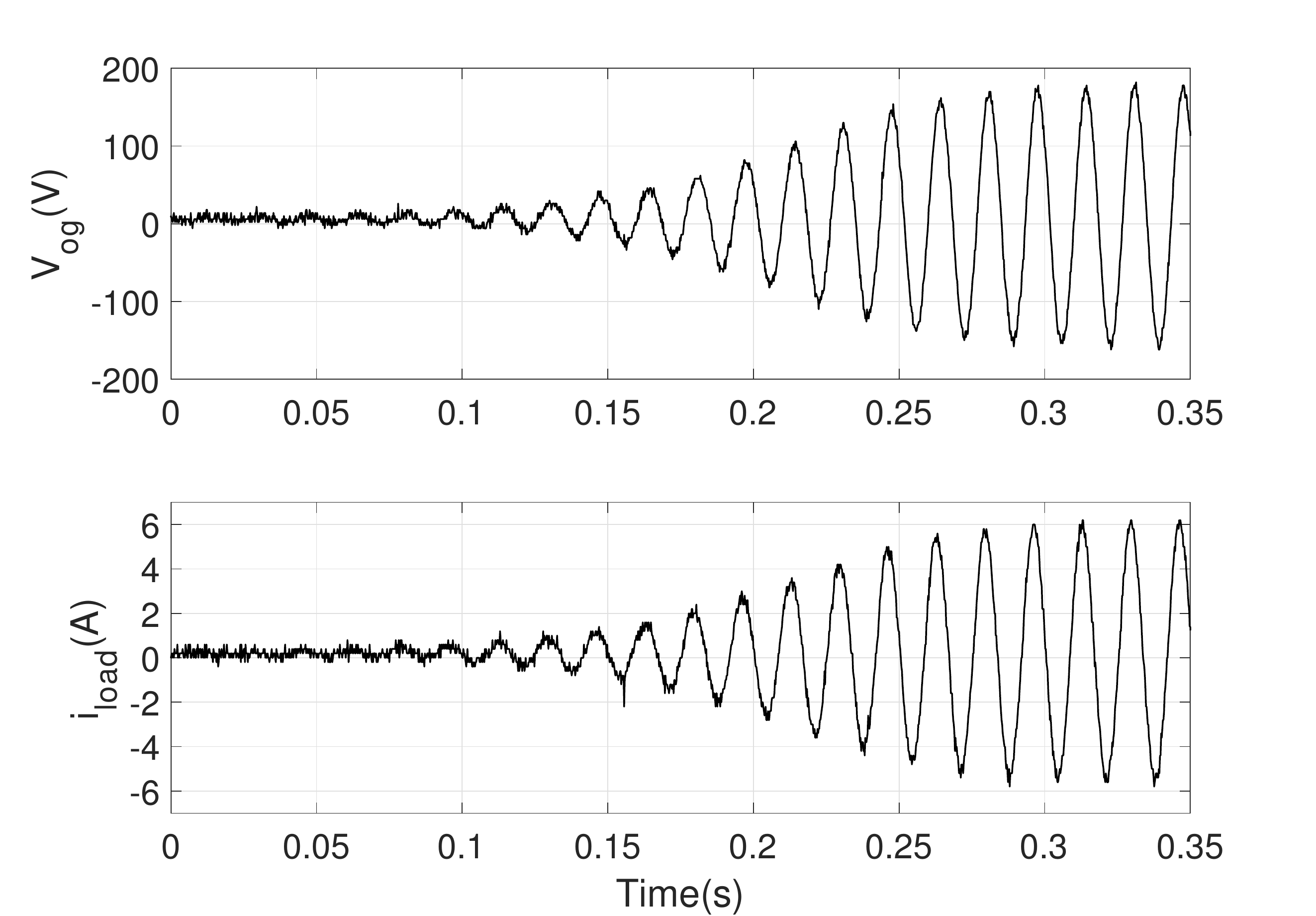}}
\caption{Black start of inverter $\#1$ under {$500$} W load. dVOC with $p^\star={500}$~W, $q^\star=-{125}$~var is used to {black-start} the grid. Since $0$ is an unstable {equilibrium} for dVOC, measurement noise is enough to drive the voltage away from $0$ to the nontrivial sinusoidal solution.}
\label{fig:ex1blackstart}
\end{figure}

\begin{figure}[t!]
\centerline{\includegraphics[width = 1.0\columnwidth]{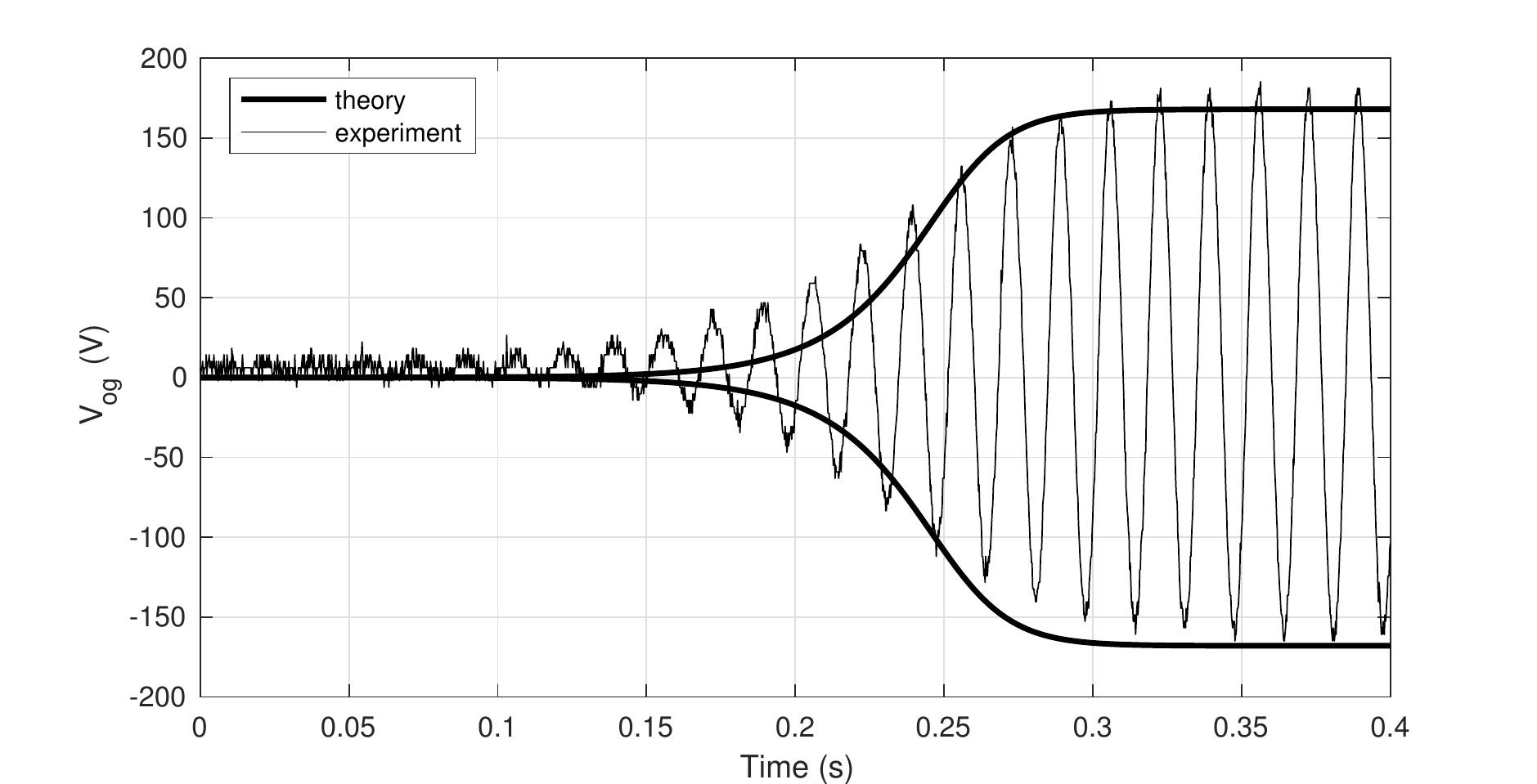}}
\caption{Comparison of the theoretical and experimental curve for the voltage magnitude during black start of inverter $\#1$ under {$500$}~W load.}
\label{fig:ex1blackstart.comp}
\end{figure}

Black start capability of grid forming inverters is a critical component for restoration after blackout to secure system resiliency in inverter-dominated power systems. Using a virtual oscillator, dVOC inverters can black-start a grid. Fig.~\ref{fig:ex1blackstart} shows that inverter \#1 initiates--i.e.,  black starts--the grid under {$500$-W} resistive load; the grid voltage is gradually established by the oscillator as expected. Due to the reactive components{,} i.e., the LCL filter connected to the load side in the system representing the backbone of an electric power system, the inverter \#1 black starts the grid under {$500$}~W and {$250$}~var (reactive components from two LCL filters) load condition. dVOC is designed such that the synchronization dynamics are much faster than the voltage dynamics. This implies that, if we consider %~\eqref{eq.droop.q.approx} 
\eqref{eq.final.droop} during black-start the second term dominates i.e.,
\[
\ddt \|v_i\| \approx  \frac{\eta \alpha}{v_i^{\star2}} \left( v_i^{\star2} -\|v_i\|^2  \right)\|v_i\|.
\]
whose solution is given by
\begin{equation}\label{black_start}
\| v_i(t) \| \approx  \frac{v_i^\star h_0 \mathrm e ^{\eta\alpha t}}{\sqrt{h_0^2 \mathrm e ^{2\eta\alpha t}+1}},
\end{equation}
where 
\[
h_0 := \frac{\|v_i(0)\|}{ \sqrt{|{\|v(0)\|^2}-v_i^{2\star}|}}.
\]
The derivation of~\eqref{black_start} is provided in the Appendix. Equation~\eqref{black_start} gives an indication of the black-start evolution of the voltage and can be used to predict the rise-time of the voltage magnitude. {Fig.~\ref{fig:ex1blackstart.comp}} shows a comparison of~\eqref{black_start} with the experimental setup during the black start of inverter $\#1$.

\begin{figure}[t!]
\centerline{\includegraphics[width = 1.0\columnwidth]{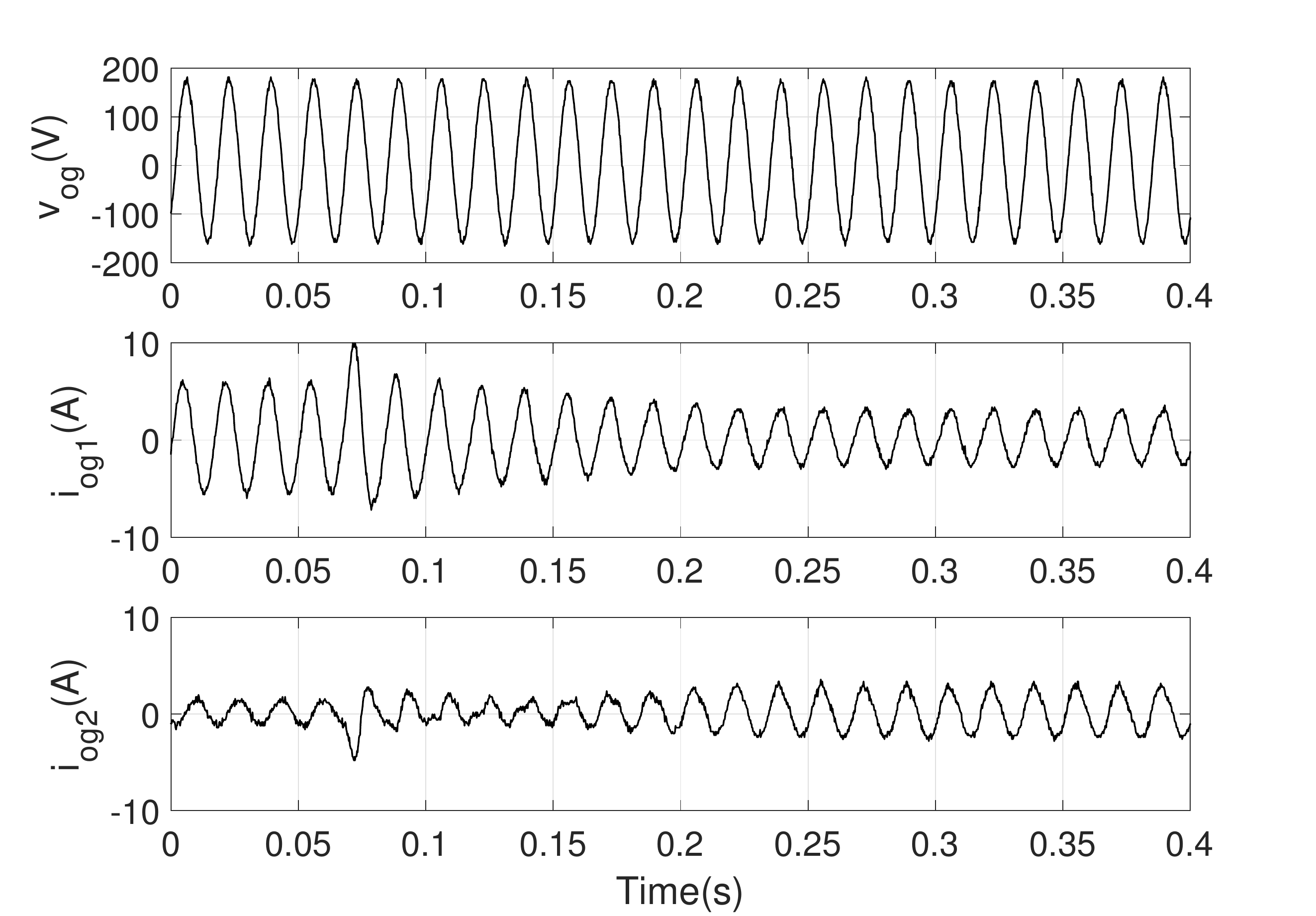}}
\caption{Connecting inverter $\#2$ at {$t=0$ sec} while inverter $\#1$ is regulating the grid under {$500$}~W load. $p^\star = {500}$~W for inverter \#1 and \#2.}
\label{fig:ex2addinginverter}
\end{figure}

\subsection{{Synchronization} of Multiple Inverters and Load Transients}

\begin{figure}[t!]
\centerline{\includegraphics[width = 1.0\columnwidth]{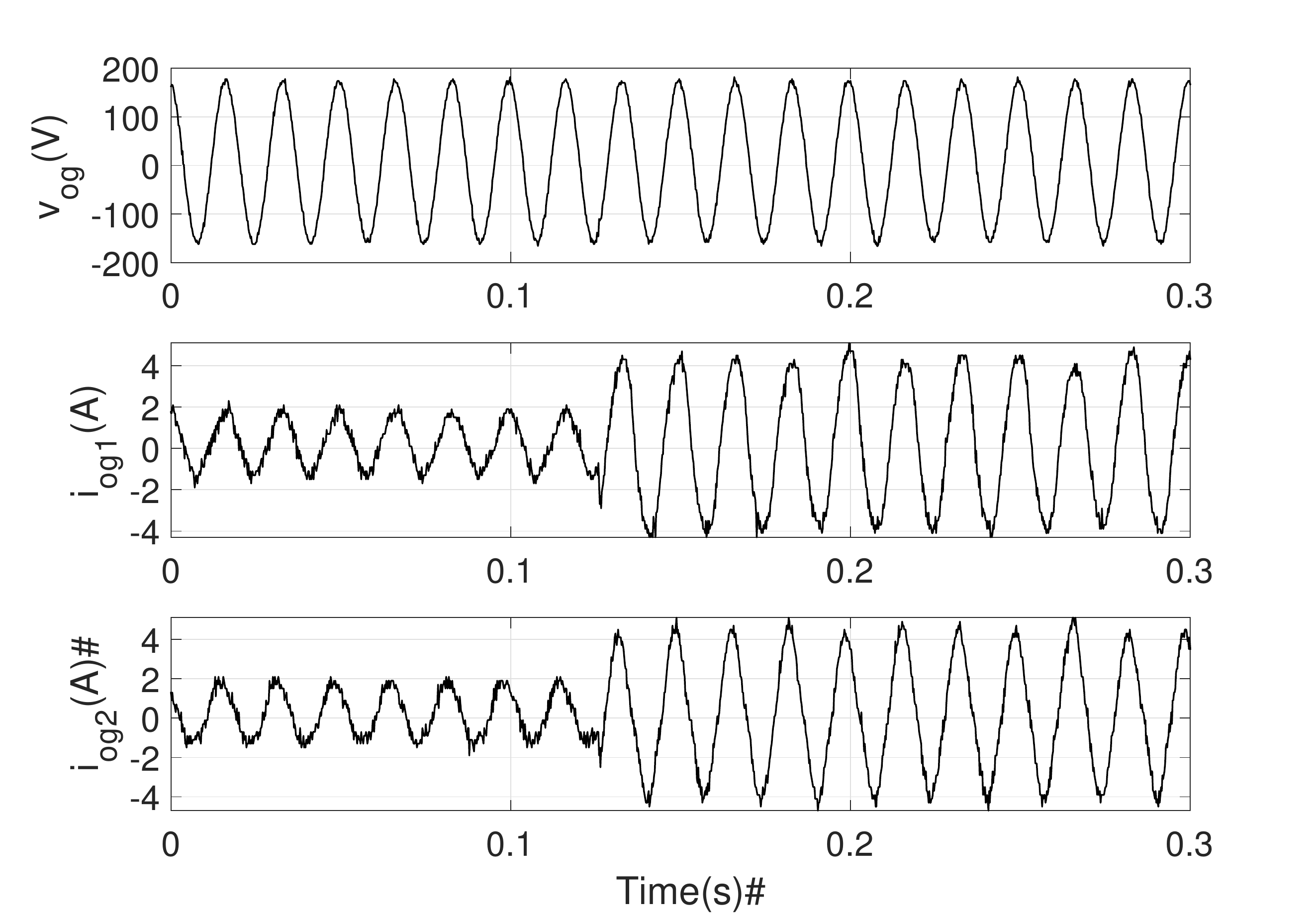}}
\caption{{$250$}~W to {$750$}~W transient with two inverters active. $p^\star$ = {$500$}~W for inverter \#1 and \#2.}
\label{fig:ex3loadtransient}
\end{figure}

As the electric grid is encountering frequent addition and subtraction of variable inverter-based renewable energy resources without significant system inertia from synchronous {machines}, the inverter controller should be capable of dynamic synchronization e.g., dynamic voltage and frequency regulation and load sharing. To verify these critical performance of grid-forming inverters, the experiment in Fig.~\ref{fig:ex2addinginverter} demonstrates a dynamic response of two grid-forming inverters; an inverter ($\#2$) is added to the grid {at $t=0$ sec} while the other $(\#1)$ is maintaining the grid. It verifies i) synchronization of multiple inverters and ii) dynamic load sharing. As the inverter $\#2$, whose LCL filter is fed by inverter $\#1$ before the transient ({i.e., the LCL filter acts as a capacitive load}), is added to the system, the two inverters synchronize within {$150$}~ms ({$10$} cycles) without significant over-current and even load sharing (same $p^\star$ set-points).
Fig.~\ref{fig:ex3loadtransient} displays operation under a load transient with two inverters active to collectively regulate the grid. Significantly different from the slow dynamic performance of conventional droop-controlled inverter operation, dVOC enables almost instantaneous dynamic load sharing without current overshoot and long settling time.

\subsection{Set-point Updates}
\begin{figure}[t!]
\centerline{\includegraphics[width = 1.0\columnwidth]{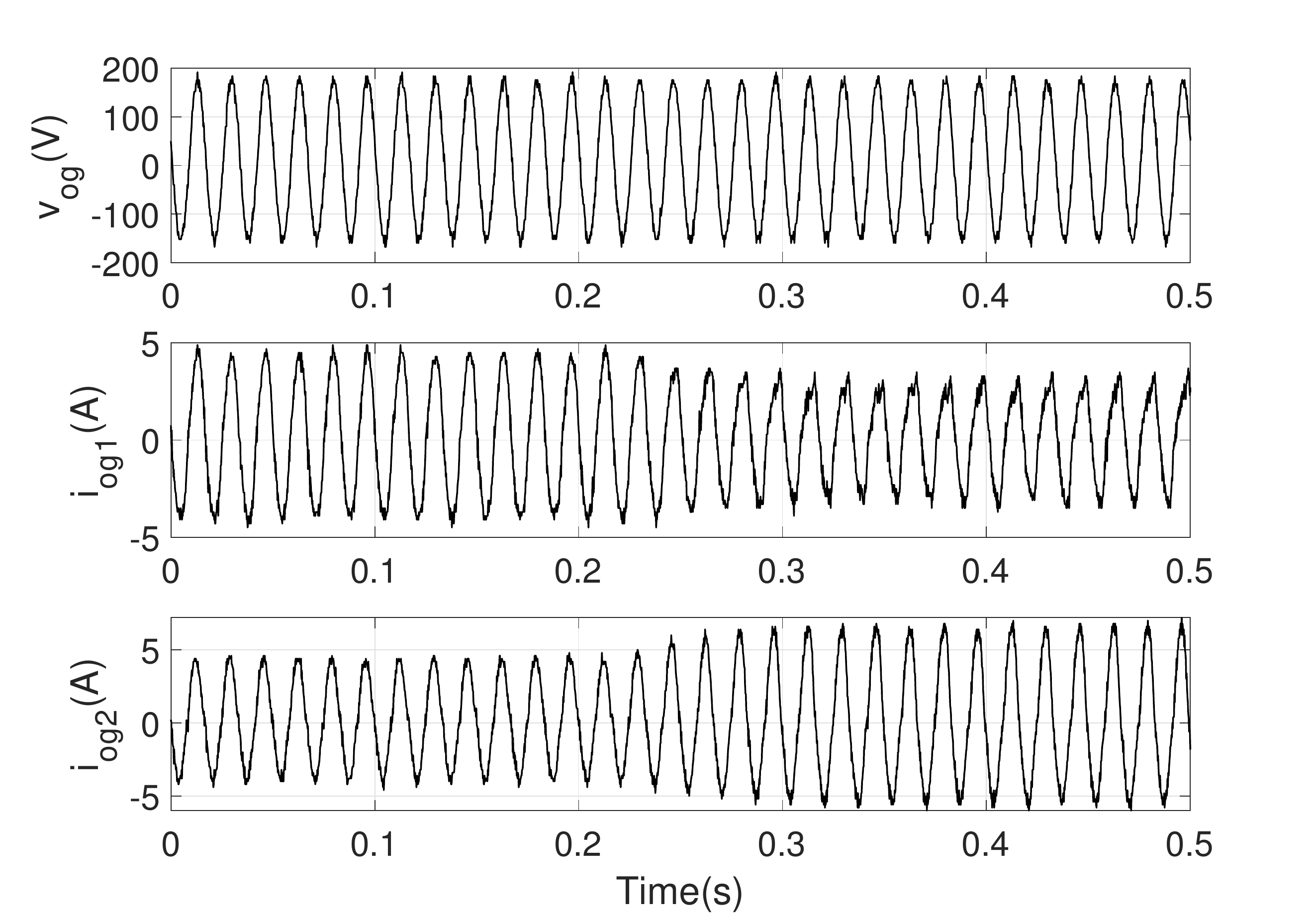}}
\caption{Change of set-point: $p_1^\star$ of inverter \#2 updated from {$250$}~W to {$500$}~W and $p_2^\star={250}$~W unchanged for inverter \#1.}
\label{fig:ex4setpointupdate}
\end{figure}

Dynamic set-point dispatch {has received} much attention due to the need for real-time power flow optimization. In an inverter-dominant grid with pervasive renewable generation, the benefit of real-time dispatch-ability will be even more striking than it is today. Because of the increased variability in generation capacity of renewable sources there will be a need of using the available resources and storage devices optimally. In Fig.~\ref{fig:ex4setpointupdate} we demonstrate a power set-point update during operation of dVOC controllers. The aim was to reconfigure the power generation profile, e.g. to optimize the power flow. As inverter \#2 changes its active power set-point $p^\star$ from {$250$}~W to {$500$}~W, the load sharing between Inverters \#1 and \# 2 is changed from {$375$}~W:{$375$}~W to {$250$}~W:{$500$}~W. In addition, the grid frequency is recovered to the nominal {$60$}~Hz since the loading condition matches to the nominal generation, verifying {that the system can be reconfigured in real time}. The ability of reconfiguring the power flow makes dVOC an extremely promising solution for an inverter-dominant system. 

\section{Conclusion and Outlook}\label{Sec:v}
%*
This paper discussed the dVOC strategy and verified its grid-forming functionality in an inverter-dominant electric power grid. Using only local information, the dVOC inverters achieve almost instantaneous dynamic synchronization and load sharing. The analysis also verified the dVOC's embedded nonlinear droop law and how the inverters can be dispatched to optimize the power flow with programmable set-points. Synchronization and dispatchability were then verified by experimental results on a custom-built hardware setup, hinting that dVOC is a promising candidate for future grid applications.

\appendix

In the following we show that~\eqref{black_start} solves the differential equation
\[
\ddt \|v_i\| =  \frac{\eta \alpha}{v_i^{\star2}} \left( v_i^{\star2} -\|v_i\|^2  \right)\|v_i\|.
\]
Making the change of variable $y = \frac{\|v_i\|}{ v_i^{\star}}$ we obtain
\[
\ddt y =  {\eta \alpha} \left( y -y^3  \right).
\]
By solving 
\[
\int_{y(0)}^{y(t)} \frac{1}{y -y^3}\mathrm d y = \eta \alpha \int_0^t \mathrm d\tau,
\]
and taking the exponent of both sides we obtain
\[
\frac{|y(t)|}{\sqrt{|y^2(t)-1|}} = \frac{|y(0)|}{\sqrt{|y^2(0)-1|}} \mathrm e ^{\eta\alpha t},
\]
by inverting the left-hand side in the interval $[0,\,1]$ and changing the variable back to $\|v_i\|$ we obtain~\eqref{black_start}.

\section*{Acknowledgements}
The authors would like to thank Dr. Miguel Rodriguez and Mr. Ramanathan Thiagarajan for their assistance in the inverter hardware development.

\bibliographystyle{ieeetr}
%\bibliography{references}
\addcontentsline{toc}{section}{\refname}
\bibliography{bib_file.bib,references.bib}

\end{document}